\documentclass[11pt]{article}

\usepackage{amsmath, amssymb, amsthm, amscd}
\usepackage{cite}

\newcommand{\ds}{\displaystyle}

\theoremstyle{plain}
\newtheorem{theorem}{Theorem}[section]
\newtheorem{lemma}[theorem]{Lemma}
\newtheorem{proposition}[theorem]{Proposition}
\newtheorem{corollary}[theorem]{Corollary}
\newtheorem*{conjecture}{Conjecture}
\newtheorem*{mainresult}{Main Result}

\theoremstyle{definition}
\newtheorem*{definition}{Definition}

\newtheorem*{remark}{Remark}

\newtheorem*{acknowledgments}{Acknowledgments}

\numberwithin{equation}{section}
\numberwithin{theorem}{section}

\begin {document}

\title{Arbitrary $p$-Gradient Values}
\author{Nathaniel Pappas\\University of Virginia}
\maketitle

\begin{abstract}
For any prime number $p$ and any positive real number $\alpha$, we construct a finitely generated group $\Gamma$ with $p$-gradient equal to $\alpha$.  This construction is used to show that there exist uncountably many pairwise non-commensurable groups that are finitely generated, infinite, torsion, non-amenable, and residually-$p$.
\end{abstract}


\footnote{This paper is to appear in J. Group Theory}

\section{Introduction}
Let $G$ be a finitely generated group and $d(G)$ denote the minimum number of generators of G. Recall the Schreier index formula: Let $H$ be a finite index subgroup of a finitely generated group $G$. Then $d(H) -1 \leq (d(G)-1)[G:H]$ and if $G$ is free of finite rank, then $H$ is free and $d(H)-1 = (d(G)-1)[G:H]$.  The rank gradient of a finitely generated group is, in a sense, a measure of how far the Schreier index formula is from being an equality rather than an inequality. Though this is an interesting question from a group-theoretic standpoint, Mark Lackenby first introduced the rank gradient as a means to study $3$-manifold groups \cite{LacExpanders}.
 
The \textit{absolute rank gradient} of $G$ is defined by 
\[RG(G) = \inf_{[G:H] < \infty} \frac{d(H) - 1}{[G : H]}\]  
where the infimum is taken over all finite index subgroups $H$ of $G$. 

It will be evident later that the rank gradient of a group is sometimes difficult both to work with and to calculate.  It is often more convenient to compute the rank gradient of the pro-$p$ completion $G_{\widehat{p}}$ of the group $G$ for some fixed prime $p$.  When dealing with profinite groups we use the notion of topologically finitely generated instead of (abstractly) finitely generated. The $p$-gradient of the group $G$, denoted $RG_{p}(G)$, can be defined as the rank gradient of $G_{\widehat{p}}$. The $p$-gradient is also referred to in the literature as the mod-$p$ rank gradient or mod-$p$ homology gradient.  A more explicit definition of $p$-gradient is provided in Section~\ref{p Gradient}.

Since Lackenby first defined rank gradient of a finitely generated group \cite{LacExpanders}, the following conjecture has remained open:

\begin{conjecture}
For every real number $\alpha > 0$ there exists a finitely generated group $\Gamma$ such that $RG(\Gamma) = \alpha$.
\end{conjecture}

The aim of this paper is to prove the analogous question for $p$-gradient:

\begin{mainresult}  For every real number $\alpha > 0$ and any prime $p$, there exists a finitely generated group $\Gamma$ such that $RG_{p}(\Gamma) = \alpha$.
\end{mainresult}

Given a prime $p$ and an $\alpha > 0 \in \mathbb{R}$, consider a free group $F$ of finite rank greater than $\alpha +1$.   Take the set of all residually-$p$ groups that are homomorphic images of $F$ that have $p$-gradient greater than or equal to $\alpha$ and partially order this set by $G_{1} \succcurlyeq G_{2}$ if $G_{1}$ surjects onto $G_{2}$.  Then by a Zorn's Lemma argument this set has a minimal element, $\Gamma$. We show $RG_{p}(\Gamma) = \alpha$ by contradiction by constructing an element which is less than $\Gamma$ with respect to the partial order.  To construct this new smaller element, Theorem~\ref{pRGlowerbound} is used, which was proved using slightly different language and a different method by Barnea and Schlage-Puchta in \cite{BarneaPuchta}, but is formulated and proved independently here as well.    

The methods used to prove this result require and are similar to those used by Schlage-Puchta in his work on $p$-deficiency and $p$-gradient \cite{Puchta} and Osin in his work on rank gradient \cite{OsinRGTorsion}.  To prove the above set has a minimal element, we will use direct limits of groups and show the relationship between the $p$-gradient of each group in the direct limit and the $p$-gradient of the limit group. This idea (Lemma~\ref{RGlimit}) was inspired by Pichot's similar result for $L^{2}$-Betti numbers \cite{Pichot}.  It is known that for a finitely generated, residually finite, infinite group the rank gradient is always greater than or equal to the $L^{2}$-Betti number, which provides a useful relationship between these two group invariants \cite{OsinRGTorsion}.

One of the primary goals of Osin's \cite{OsinRGTorsion} and Schlage-Puchta's \cite{Puchta} papers was to provide a simple construction of non-amenable, torsion, residually finite groups.  The construction given in this paper shows that there exist such groups with arbitrary $p$-gradient (Theorem~\ref{arbitrarytorsionRGp}). A simple consequence of this result is that there exist uncountably many pairwise non-commensurable groups that are finitely generated, infinite, torsion, non-amenable, and residually-$p$.  The fact that the groups are non-commensurable uses the $p$-gradient and is almost immediate from the construction, which shows another way in which the $p$-gradient can be a useful tool.  

\begin{acknowledgments}
The author would like to thank his advisor Mikhail Ershov for his help with the present material and earlier drafts of this paper.  The author would also like to thank the anonymous referee for pointing out a minor mathematical issue and helpful comments which improved the exposition.
\end{acknowledgments}

\section{Rank Gradient and $p$-Gradient}
\label{p Gradient}
In this section, some useful results concerning rank gradient and $p$-gradient are collected, which will be used to prove the main result.  

\begin{theorem}
\label{RGFiniteIndex}
Let $G$ be a finitely generated group and let $H$ be a finite index subgroup.  Then $RG(G) = \frac{RG(H)}{[G:H]}$.  If $G$ is finite, then $RG(G) = -\frac{1}{|G|}$.
\end{theorem}

\begin{proof}
By the Schreier index formula, $\frac{d(K)-1}{[G:K]}\geq \frac{d(L)-1}{[G:L]}$ for any subgroups $L \leq K$ of $G$  such that $L$ is finite index in $G$.  This fact is used in the first equality given below.

Fix a finite index subgroup $H \leq G$.  Since any finite index subgroup $K$ of $G$ contains the finite index subgroup $K \cap H$ of $H$ and any finite index subgroup $L$ of $H$ is $L = K \cap H$ for some finite index subgroup $K$ of $G$, then 
$\ds RG(G) = \inf_{[G:K] < \infty} \frac{d(K)-1}{[G:K]} = \inf_{[G:K] < \infty} \frac{d(K \cap H)-1}{[G:K \cap H]} = \inf_{[G:K] < \infty} \frac{d(K \cap H)-1}{[G:H][H:K \cap H]} = \frac{1}{[G:H]} \inf_{[H:L] < \infty} \frac{d(L)-1}{[H:L]} = \frac{RG(H)}{[G:H]}.$  

If $G$ is finite, then using $H = \{1\}$ implies $RG(G) = \frac{-1}{|G|}$. 
\end{proof}

As the following proposition shows, it is not difficult to produce groups with rational rank gradient.  Whether an irrational number can be the rank gradient of some finitely generated group remains an open question. We will show later that for every prime $p$, every positive real number is the $p$-gradient for some finitely generated group.

\begin{proposition}
Let $q >0 \in \mathbb{Q}$.  There exists a finitely presented group $G$ such that $RG(G) = q$.
\end{proposition}

\begin{proof}
Write $q = \frac{m}{n}$.  Let $F_{m+1}$ be a non-abelian free group of rank $m+1$ and let $A$ be any group of order $n$.  Consider $G = F_{m+1} \times A$.  Let $\varphi:G \to A$ be the projection onto the second component and let $H = \ker{\varphi}$.  Then $F_{m+1} \simeq H$ and $[G:H] = n$.  By the Schreier index formula for free groups $RG(H) = m$.  Therefore, $RG(G) = \frac{RG(H)}{[G:H]} = \frac{m}{n}$ by Theorem~\ref{RGFiniteIndex}.
\end{proof}

As stated earlier, we can define $RG_{p}(G) = RG(G_{\widehat{p}})$. However, a more explicit definition of the $p$-gradient can be stated.

\begin{definition} 
Let $p$ be a prime. The \textit{$p$-gradient} (also called mod-$p$ homology gradient) of $G$ is defined by 
\[RG_{p}(G) = \inf \frac{d_{p}(H)-1}{[G : H]}\]
where $d_p(G) = d\left(G / [G,G]G^{p} \right)$ and the infimum is taken over all normal subgroups $H$ such that $[G : H] = p^{k}$ for some $k \in \mathbb{Z}_{\geq0}$.
\end{definition}

We will prove that a group and its pro-$p$ completion have the same $p$-gradient, which will then be used to show the $p$-gradient of a group equals the rank gradient of its pro-$p$ completion. To get this result, some facts about profinite groups must be presented.  

Let $G$ be a finitely generated group. The pro-$p$ completion of $G$ for some prime $p$ will be denoted by $G_{\widehat{p}}$. Let $d(G)$ denote the minimal number of abstract generators of a group $G$ if the group is not profinite and the minimal number of topological generators if the group is profinite.  If a group is profinite, the term ``finitely generated" will be used to mean topologically finitely generated.  The reader is referred to any standard text in profinite groups for the basic results used in this section \cite{DDMS}, \cite{wilson}.  

When dealing with pro-$p$ completions of a group, it is often convenient to assume that the group is residually-$p$ since the group will imbed in its pro-$p$ completion.  To show why this type of assumption will not influence any result about the $p$-gradient, the following lemma is given.  

\begin{definition}
Let $G$ be a group and $p$ a prime.  Let $N$, the $p$-residual of $G$, be the intersection of all normal subgroups of $p$-power index in $G$.  The \textit{$p$-residualization} of $G$ is the quotient $G/N$. Note that the $p$-residualization of $G$ is isomorphic to the image of $G$ in its pro-$p$ completion $G_{\widehat{p}}$ and is residually-$p$.
\end{definition}

\begin{lemma}
\label{residppropcomp}
Let $G$ be a group and $p$ a prime number.  Let $\widetilde{G}$ be the $p$-residualization of G.  Then $RG_{p}(G) = RG_{p}(\widetilde{G})$ and $G_{\widehat{p}} \simeq \widetilde{G}_{\widehat{p}}$.
\end{lemma}

\begin{proof}
There is a bijective correspondence between normal subgroups of $p$-power index in $\widetilde{G}$ and normal subgroups of $p$-power index in $G$.  Let the correspondence be $\widetilde{H} \leftrightarrows H$ with $\widetilde{H} \leq \widetilde{G}$ and $H \leq G$.  Then it is easy to show that $[\widetilde{G}:\widetilde{H}] = [G:H]$ and $d_{p}(\widetilde{H}) = d_{p}(H)$.  Thus,  $RG_{p}(\widetilde{G}) = RG_{p}(G)$.  

By the inverse limit definition of pro-$p$ completion and the fact that $\widetilde{G}/\widetilde{H} \simeq (G/N) / (H/N) \simeq G/H$ it follows that $G_{\widehat{p}} \simeq \widetilde{G}_{\widehat{p}}$.
\end{proof}

The following fact is well known.

\begin{proposition}
\label{subnormal base for pro-p}
Let $G$ be a group and $p$ a prime number.  The set of subnormal subgroups of $p$-power index form a base of neighborhoods of the identity for the pro-$p$ topology on $G$. 
\end{proposition}

With the following proposition, we will be able to prove that a group and its pro-$p$ completion have the same $p$-gradient.  Parts of this proposition can be found in an exercise in \cite{DDMS}. 

\begin{proposition}
\label{bijectivecorespondence}
Let $G$ be a finitely generated group and $p$ a prime.  Let $\varphi: G \to G_{\widehat{p}}$ be the natural map from $G$ to its pro-$p$ completion.  Let $H$ be a normal subgroup of $p$-power index of $G$.  The following hold:
\begin{enumerate}
\item $\varphi(H) = \varphi(G) \cap  \overline{\varphi(H)}$.
\item $\overline{\varphi}:G/H \to G_{\widehat{p}} / \overline{\varphi(H)}$ given by $\overline{\varphi}(xH) = \varphi(x) \overline{\varphi(H)}$ is an isomorphism.
\item There exists an index preserving bijection between normal subgroups of $p$-power index in $G$ and open normal subgroups of $G_{\widehat{p}}$.  
\item  $\overline{\varphi(H)} \simeq H_{\widehat{p}}$ as pro-$p$ groups.
\item $\ds RG(G_{\widehat{p}}) = \frac{RG(H_{\widehat{p}})}{[G:H]}$.
\end{enumerate}  
\end{proposition}

\begin{proof}
Parts (1)-(3) are proved in Proposition 3.2.2 of Ribes and Zalesskii \cite{Ribes}.   
\item{4)}  For notational simplicity, assume $G$ is residually-$p$ and thus $\varphi$ is injective.  The case of $G$ not residually-$p$ is proved similarly. 
It is only necessary to show that the pro-$p$ topology on $G$ induces the pro-$p$ topology on the subspace $H$ of $G$.  By Proposition~\ref{subnormal base for pro-p}, subnormal subgroups of $p$-power index in $G$ form a base for the pro-$p$ topology on $G$.  If $K$ is subnormal of $p$-power index in $H$ it implies that $K$ is subnormal of $p$-power index in $G$.  This implies that the subspace topology and the pro-$p$ topology on $H$ are the same.  Therefore, $\overline{H} \simeq H_{\widehat{p}}$ as pro-$p$ groups.   
\item{5)} By (2) and (4) we know $G/H \simeq G_{\widehat{p}}/ H_{\widehat{p}}$ and therefore $[G:H] = [G_{\widehat{p}}:H_{\widehat{p}}]$.  Thus, by Theorem~\ref{RGFiniteIndex} $RG(G_{\widehat{p}}) = \frac{RG(H_{\widehat{p}})}{[G_{\widehat{p}}:H_{\widehat{p}}]} = \frac{RG(H_{\widehat{p}})}{[G:H]}$.
\end{proof}

\begin{theorem}
\label{propcompRGp}
If $G$ is a (topologically) finitely generated pro-$p$ group, then $RG_{p}(G) = RG(G)$.
\end{theorem}

\begin{proof}
In a finitely generated pro-$p$ group all finite index normal subgroups are open normal subgroups and have index a power of $p$ \cite{DDMS}.  Moreover, if $H$ is a finite index subgroup of $G$, then $H$ is also a finitely generated pro-$p$ group.  The Frattini subgroup of a finitely generated pro-$p$ group $H$ is $\Phi(H) = [H,H]H^{p}$ and by standard facts about finitely generated pro-$p$ groups, $d_{p} (H) = d(H / \Phi(H)) = d(H).$  
Therefore, 
\[RG_{p}(G) = \inf_{\substack{{H \trianglelefteq G} \\ [G:H] = p^{k}}} \frac{d_{p}(H)-1}{[G:H]} 
 =  \inf_{[G:H] \leq \infty} \frac{d(H)-1}{[G:H]} = RG(G).\]
 \end{proof}

It is now possible to prove the relationship between the $p$-gradient of a group and its pro-$p$ completion.  

\begin{theorem}
\label{pRGRelationTheorem}
Let $G$ be a finitely generated group and $p$ a fixed prime.  Let $G_{\widehat{p}}$ be the pro-$p$ completion of $G$. Then $RG_{p}(G) = RG_{p}(G_{\widehat{p}}) = RG(G_{\widehat{p}})$. 
\end{theorem}

\begin{proof}
Assume that $G$ is residually-$p$.  By Proposition~\ref{bijectivecorespondence} the natural injective map $\varphi: G \to G_{\widehat{p}}$ induces an index preserving bijection $H \to \overline{H} \simeq H_{\widehat{p}}$ between the normal subgroups of $p$-power index in $G$ and the normal subgroups of $p$-power index in $G_{\widehat{p}}$. Proposition~\ref{bijectivecorespondence} also implies that $d_{p}(\overline{H}) = d_{p}(H_{\widehat{p}})$ for all $p$-power index normal subgroups $H \trianglelefteq G$. It is not difficult to show that $d_{p}(K) = d_{p}(K_{\widehat{p}})$ holds for any finitely generated group $K$.  Therefore, $RG_{p}(G) = RG_{p}(G_{\widehat{p}})$. If $G$ is not residually-$p$, let $\widetilde{G}$ be the $p$-residualization of G.  By Lemma~\ref{residppropcomp}, $RG_{p}(G)=RG_{p}(\widetilde{G})$ and $\widetilde{G}_{\widehat{p}} \simeq G_{\widehat{p}}$.  The first equality follows.

The fact that $RG_{p}(G) = RG(G_{\widehat{p}})$, where $G_{\widehat{p}}$ is the pro-$p$ completion of $G$, follows by the above remarks and Theorem~\ref{propcompRGp}. 
\end{proof}

\begin{remark}  
Nikolov and Segal proved Serre's conjecture on finitely generated profinite groups.  That is, in a finitely generated profinite group all finite index subgroups are open \cite{NikolovSegal}.
\end{remark}

The above two theorems provide some useful corollaries.  

\begin{corollary}
\label{finitepRG}
If $G$ is a finite group, then $RG_{p}(G) = -\frac{1}{|G_{\widehat{p}}|}$.
\end{corollary}

\begin{proof}
If $G$ is finite, then so is $G_{\widehat{p}}$ and thus $RG_{p}(G) = RG_{p}(G_{\widehat{p}}) = RG(G_{\widehat{p}}) = -\frac{1}{|G_{\widehat{p}}|}$ by Theorem~\ref{RGFiniteIndex}.
\end{proof}

\begin{theorem}
\label{pindexRGp}
Fix a prime $p$ and let $G$ be a finitely generated group.  Assume $H \leq G$ is a $p$-power index subnormal subgroup.  Then $RG_{p}(G) = \frac{RG_{p}(H)}{[G:H]}$.
\end{theorem}

\begin{proof}
If $H$ in normal in $G$ we are done by Proposition~\ref{bijectivecorespondence}.5.  If $H$ is not normal, then induct on the subnormal length.    
\end{proof}

\section{Groups With Arbitrary $p$-Gradient Values}
In this section we will prove the main result, that is, we construct a finitely generated group $\Gamma$ with $RG_{p}(\Gamma)=\alpha$ for each $\alpha>0 \in \mathbb{R}$.  To prove this, we need some technical results.  

The following lemma is similar to Lemma 2.3 of Osin \cite{OsinRGTorsion} concerning deficiency of a finitely presented group.

\begin{lemma}
\label{AddingRelationLemma}
Let $G$ be a finitely generated group and fix a prime $p$.  Let $x$ be some non-trivial element of $G$.  Let $H$ be a finite index normal subgroup of $G$ such that $x^{m} \in H$, but no smaller power of $x$ is in $H$.  Let $\pi:G \to G / \langle x^{m} \rangle^{G}$ be the standard projection homomorphism.
\begin{enumerate}
\item If $T$ is a right transversal for $\langle x \rangle H$ in $G$, then $\langle x^{m} \rangle^{G} = \\ \langle tx^{m}t^{-1} \mid t \in T \rangle^{H}$.
\item If $H=\langle Y \mid R \rangle$, then $\pi(H) = \langle Y \mid R \cup \{tx^{m}t^{-1} \mid t \in T\} \rangle$.  
\item $\ds |T| = \frac{[G:H]}{m}$.
\item If $\ds \mathfrak{q}(H) = \frac{d_{p}(H)}{[G:H]}$, then $\mathfrak{q}(\pi(H)) \geq \mathfrak{q}(H) - \frac{1}{m}$.   
\end{enumerate}
\end{lemma}
\vspace{5 pt}

\begin{proof}
Since $x^{m}$ is in $H$, then $[\pi(G):\pi(H)] = [G:H].$
\begin{enumerate}
\item  This is a standard computation.

\item  This holds by (1) and the fact that $\pi(H) = H / (H \cap \langle x^{m} \rangle^{G}) = H / \langle x^{m} \rangle^{G}$, since $x^{m} \in H$ and $H$ is normal in $G$.

\item Since $H \subseteq \langle x \rangle H \subseteq G$, then $[G:H] = [G: \langle x \rangle H][\langle x \rangle H:H]$.  Therefore, $ |T| = [G:\langle x \rangle H] = \frac{[G : H ]}{[\langle x \rangle H : H ]}$.  Since $x^{m} \in H$ but no smaller power of $x$ is in $H$, then $V=\{1, x, x^{2}, \dots, x^{m-1}\}$ is a transversal for $H$ in $\langle x \rangle H$ and thus $[\langle x \rangle H:H] = m$.  Therefore, $|T| = \frac{[G:H]}{m}$. 

\item  First, note that (2) and (3) imply that a presentation for $\pi(H)$ is obtained from a presentation for $H$ by adding in $\frac{[G:H]}{m}$ relations.  Now, $\mathfrak{q}(\pi(H)) \geq \mathfrak{q}(H) - \frac{1}{m}$ if and only if $d_{p}(\pi(H)) \geq d_{p}(H) - \frac{[G:H]}{m}.$ 
If $H$ has presentation $H = \langle Y \mid R \rangle$ then $\pi(H)$ has presentation \\
$\pi(H) =  \langle Y \mid R \cup \{tx^{m}t^{-1} \hspace{4 pt} \text{for all} \hspace{4 pt} t \in T\} \rangle.$
For notational simplicity let $C = \{[y_{1},y_{2}] \mid y_{1}, y_{2} \in Y\}$.  Then, 
 \[H / ([H,H]H^{p}) = \langle Y \mid R, \hspace{4 pt} C, \hspace{4 pt} w^{p} \hspace{4 pt} \text{for all} \hspace{4 pt} w \in F(Y) \rangle\] 
 where $F(Y)$ is the free group on $Y$ and 
 \begin{align*}
 \pi(H) / ([\pi(H),\pi(H)]\pi(H)^{p}) = \langle Y \mid & \hspace{4 pt} R, \hspace{4 pt} C, \hspace{4 pt} w^{p} \hspace{4 pt} \text{for all} \hspace{4 pt} w \in F(Y), \\
 &tx^{m}t^{-1} \text{for all} \hspace{4 pt} t \in T \rangle.
 \end{align*}
 
Therefore, a presentation for $\pi(H) / ([\pi(H),\pi(H)]\pi(H)^{p})$ is obtained from a presentation for $H / ([H,H]H^{p})$ by adding in $\frac{[G:H]}{m}$ relations.

 \textbf{Note:} For any group G, $G/([G,G]G^{p})$ can be considered as a vector space over $\mathbb{F}_{p}$ and therefore $d_{p}(G)$ is the dimension of this vector space.
 
Therefore, $\pi(H) / ([\pi(H),\pi(H)]\pi(H)^{p})$ is a vector space satisfying $\frac{[G:H]}{m}$ more equations than the vector space $H / ([H,H]H^{p})$.
Thus $d_{p}(\pi(H)) \geq d_{p}(H) - \frac{[G:H]}{m}.$
\qedhere
\end{enumerate}
\end{proof}

A lower bound for the $p$-gradient when taking the quotient by the normal subgroup generated by an element raised to a $p$-power follows by the above lemma.   

\begin{theorem}
\label{pRGlowerbound}
Let $G$ be a finitely generated group, $p$ some fixed prime, and $x \in G$.  Then $RG_{p}(G/\langle\langle x^{p^{k}} \rangle\rangle) \geq RG_{p}(G) - \frac{1}{p^{k}}.$  
\end{theorem}

\begin{proof}
\underline{Case 1:} There exists a normal subgroup $H_{0}$ of $p$-power index such that the order of $x$ in $G/H_{0}$ is at least $p^{k}$.

Since $H_{0}$ is a normal subgroup of $p$-power index, then without loss of generality we may assume that the order of $x$ in $G/H_{0}$ is exactly $p^{k}$.  Let $\overline{H}$ be a normal subgroup of $p$-power index in $\overline{G} =  G/ \langle\langle x^{p^{k}} \rangle\rangle.$  Let $H \leq G$ be the full preimage of $\overline{H}$.  Then $H$ is a $p$-power index normal subgroup in $G$ which contains $\langle\langle x^{p^{k}} \rangle\rangle$.  Let $L_{H} = H \cap H_{0}$.  Then $L_{H}$ is a normal subgroup in $G$ such that $x^{p^{k}} \in L_{H}$, $L_{H} \subseteq H$, and the order of $x$ in $G/L_{H}$ is $p^{k}$.  Note that $L_{H}$ is normal and of $p$-power index in $G$ since both $H$ and $H_{0}$ are normal and of $p$-power index.  Thus by Lemma~\ref{AddingRelationLemma}, $\mathfrak{q}(\overline{H}) \geq \mathfrak{q}(\overline{L_{H}}) \geq \mathfrak{q}(L_{H}) - \frac{1}{p^{k}},$ which by definition is greater than or equal to $RG_{p}(G) - \frac{1}{p^{k}}$.  Therefore, $\mathfrak{q}(\overline{H}) \geq RG_{p}(G) - \frac{1}{p^{k}}$.  Thus $RG_{p}(G/\langle\langle x^{p^{k}} \rangle\rangle) \geq RG_{p}(G) - \frac{1}{p^{k}}.$  

\underline{Case 2:} For every normal subgroup $H$ of $p$-power index, the order of $x$ in  $G/H$ is less than $p^{k}$.

It will be shown that $RG_{p}(G/\langle\langle x^{p^{k}} \rangle\rangle) = RG_{p}(G)$ in this case. There exists an $\ell < k$ such that $x^{p^{\ell}} \in H$ for every normal subgroup $H$ of $p$-power index in $G$.  Then $x^{p^{\ell}}$ is in the kernel of natural map from $G$ to its pro-$p$ completion $\varphi: G \to G_{\widehat{p}}$. Therefore, $x^{p^{k}} = (x^{p^{\ell}})^{p^{k-\ell}} \in \ker \varphi.$  Let $M = \langle\langle x^{p^{k}} \rangle\rangle$.  Then $M \subseteq \ker \varphi$. This implies that there is a bijective correspondence between all normal subgroup of $p$-power index in $G$ and $G/M$ given by $N \to N/M$.  Since $G/N \simeq (G/M)/(N/M)$ for all such $N$, then by the inverse limit definition of pro-$p$ completions $G_{\widehat{p}} \simeq (G/M)_{\widehat{p}}$ as pro-$p$ groups.  Therefore, $RG_{p}(G/\langle\langle x^{p^{k}} \rangle\rangle) = RG_{p}(G/M) = RG_{p}((G/M)_{\widehat{p}}) = RG_{p}(G_{\widehat{p}}) = RG_{p}(G).$
\end{proof}

\begin{remark} 
The above theorem was independently stated and proved using slightly different language and a different method by Barnea and Schlage-Puchta (Theorem 3 in \cite{BarneaPuchta}).
\end{remark}

\begin{corollary}
Let $G$ be a finitely generated group, $p$ a fixed prime, and let $x \in G$.  Then $RG_{p}(G/\langle\langle x \rangle\rangle) \geq RG_{p}(G) - 1.$  
\end{corollary}

\subsection{$p$-Gradient and Direct Limits}
Let $(I, \leq)$ be a totally ordered set with smallest element $0$ and let $\{G_{i} \mid \pi_{ij}\}$ be a direct system of finitely generated groups with  surjective homomorphisms $\pi_{ij}:G_{i} \to G_{j}$ for every $j \geq i \in I$.  

Let $G_{\infty} = \varinjlim G_{i}$ be the direct limit of this direct system.   Let $\pi_{i}:G_{i} \to G_{\infty}$ be the map obtained from the direct limit.  Because all the maps in the direct system are surjective, then so are the $\pi_{i}.$  Let $G = G_{0}$.

Another direct system $\{M_{i} \mid \mu_{ij}\}$ can be defined over the same indexing set $I$, where $M_{i}=G$ for each $i$ and $\mu_{ij}$ is the identity map.  The direct limit of this set is clearly $G = \varinjlim M_{i}$ and the map obtained from the direct limit $\mu_{i}:M_{i} \to G_{\infty}$ is the identity map.

A homomorphism $\Phi:\{M_{i} \mid \mu_{ij}\} \to \{G_{i} \mid \pi_{ij}\}$ is by definition a family of group homomorphisms $\varphi_{i}:M_{i} \to G_{i}$ such that $\varphi_{j} \circ \mu_{ij} =\pi_{ij} \circ \varphi_{i}$ whenever $i \leq j$.  Then $\Phi$ defines a unique homomorphism $\varphi = \varinjlim \varphi_{i}: \varinjlim M_{i} \to \varinjlim G_{i}$ such that $\varphi \circ \mu_{i}=\pi_{i} \circ \varphi_{i}$ fror all $i \in I$ \cite{AtiyahMacDonald}.

The surjection $\varphi_{i}:G \to G_{i}$ is the map $\pi_{0i}$ in this case. It is clear that $\varphi = \varinjlim \varphi_{i}$.  Since each $\varphi_{i}$ is surjective, it implies that $\ker{\varphi_{i}} \subseteq \ker{\varphi_{j}}$ for every $j \geq i$. In this situation,
 \[\ker{\varphi} = \varinjlim \ker{\varphi_{i}} = \bigcup_{i \in I} \ker{\varphi_{i}}.\]

Let $H \leq G$ be a subgroup.  For every $i$, let $H_{i} = \varphi_{i}(H)$.

\begin{lemma}
\label{subgrouplimits}
Keep the notation defined above.  Fix a prime $p$.  For each $K \trianglelefteq G_{\infty}$ of $p$-power index, there exists an $H' \trianglelefteq G$ of $p$-power index such that: 
\begin{enumerate}
\item $K = \varinjlim H'_{i}$.
\item $\ds [G_{\infty}:K] = \lim_{i \in I} [G_{i}:H'_{i}]$.
\item $\ds d_{p}(K) = \lim_{i \in I} d_{p}(H'_{i})$.
\end{enumerate}
\end{lemma}

\begin{proof}
Let $K \trianglelefteq G_{\infty}$ be a $p$-power index normal subgroup. Since $\varphi:G \to G_{\infty}$ is surjective then $G_{\infty} \simeq G / \ker{\varphi}$.   Let $H' = \varphi^{-1}(K)$.  Then $H'$ is normal in $G$ and since $K \simeq H' / \ker{\varphi}$ then $[G_{\infty}:K] = [G:H']$ and so $H'$ is of $p$-power index.   
\begin{enumerate}
\item $K = \varphi(H') = \varinjlim \varphi_{i}(H')  = \varinjlim H'_{i}$. 

\item  Since each $\varphi_{i}:G \to G_{i}$ is surjective, $G_{i} \simeq G / \ker{\varphi_{i}}$ and since $H'$ contains $\ker{\varphi}$, then $H'$ contains $\ker{\varphi_{i}}$ for each $i$.  Thus, $H'_{i} \simeq H' / \ker{\varphi_{i}}$.  Therefore for every $i$, 
\[G_{i} / H'_{i} \simeq G / H' \simeq G_{\infty} / K.\]
Thus, $[G_{\infty}: K] = [G_{i} : H'_{i}]$ for every $i$.

\item For any group $A$, let $Q(A) = A / [A,A]A^{p}$.  It is known that $K \simeq H' / \ker{\varphi}$ and $H'_{i} \simeq H' / \ker{\varphi_{i}}$ and therefore,
\[Q(K) \simeq H' / [H',H'](H')^{p} \ker{\varphi} \simeq Q(H') / M\]
where $M = [H',H'](H')^{p} \ker{\varphi} / [H',H'](H')^{p}$, and 
\[Q(H'_{i}) \simeq H' / [H',H'](H')^{p} \ker{\varphi_{i}} \simeq Q(H') / M_{i}\]
where $M_{i} = [H',H'](H')^{p} \ker{\varphi_{i}} /[H',H'](H')^{p}$.  Since $\ker{\varphi_{i}} \subseteq \ker{\varphi_{j}}$ for each $j \geq i$ then $M_{i} \subseteq M_{j}$ for each $j \geq i$.  Now, $Q(H')$ is finitely generated abelian and torsion and therefore is finite.  Thus $Q(H')$ can only have finitely many non-isomorphic subgroups. Since $\{M_{i}\}$  is an ascending set of subgroups, there must exist an $n \in I$ such that $M_{i} = M_{n}$ for every $i \geq n$.  Since $\ker{\varphi_{i}} \subseteq \ker{\varphi_{j}}$ for each $j \geq i$ and $\bigcup \ker{\varphi_{i}} = \ker{\varphi}$, we know that $M_{i} \subseteq M_{j}$ for every $j \geq i$ and $\bigcup M_{i} = M$.   Therefore, $M = \bigcup M_{i} = M_{n}$.  Thus for each $i \geq n$,  $M = M_{i}$.  

Therefore, $Q(K) \simeq Q(H'_{i})$ for each $i \geq n$ which implies $d_{p}(K) = d_{p}(H'_{i})$ for each $i \geq n$.  Thus, $\ds d_{p}(K) = \lim_{i \in I} d_{p}(H'_{i})$.
\qedhere
\end{enumerate}
\end{proof}

The following lemma is similar to Pichot's related result for L$^{2}$-Betti numbers where convergence is in the space of marked groups \cite{Pichot}. 

\begin{lemma}
\label{RGlimit}
For each prime $p$, $\limsup RG_{p}(G_{i}) \leq RG_{p}(G_{\infty})$.  
\end{lemma}

\begin{proof}  Fix a prime $p$.  Let $K \trianglelefteq G_{\infty}$ be a normal subgroup of $p$-power index.  By Lemma~\ref{subgrouplimits}  we obtain the subgroups $H'$ and $H'_{i}$ for each $i$. Now, 
\[ \limsup RG_{p}(G_{i}) = \limsup \inf_{\substack{N \trianglelefteq G_{i} \\ \text{$p$-power}}} \frac{d_{p}(N)-1}{[G_{i}:N]} \leq \limsup \frac{d_{p}(H'_{i})-1}{[G_{i}:H'_{i}]}\] 
and by Lemma~\ref{subgrouplimits} 
\[ \limsup \frac{d_{p}(H'_{i})-1}{[G_{i}:H'_{i}]} = \lim_{i \in I} \frac{d_{p}(H'_{i})-1}{[G_{i}:H'_{i}]} = \frac{d_{p}(K)-1}{[G_{\infty}:K]}.\]
Therefore, for each $K \trianglelefteq G_{\infty}$ of $p$-power index, $\limsup RG_{p}(G_{i}) \leq \frac{d_{p}(K)-1}{[G_{\infty}:K]}$.  This implies $\limsup RG_{p}(G_{i}) \leq RG_{p}(G_{\infty})$.
\end{proof}

\subsection{The Main Result}
It is now possible to prove the main result that every nonnegative real number is realized as the $p$-gradient of some finitely generated group.  

\begin{theorem}
\label{arbitraryRGp}
\textbf{(Main Result)}
For every real number $\alpha >0$ and any prime $p$, there exists a finitely generated group $\Gamma$ such that $RG_{p}(\Gamma) = \alpha$.
\end{theorem}

\begin{proof}
Fix a prime $p$ and a real number $\alpha>0$.  Let $F$ be the free group on $\lceil \alpha \rceil +1$ generators.  Let 
\[ \Lambda = \{G \mid \text{$F$ surjects onto $G$, $G$ is residually-$p$, and $RG_{p}(G)$} \geq \alpha\}.\]
Since for any free group $d(F) = d_{p}(F)$ it is clear that $RG_{p}(F) = $ rank$(F)-1$ and therefore, $\Lambda$ is not empty since $F$ is in $\Lambda$.  $\Lambda$ can be partially ordered by $G_{1} \succcurlyeq G_{2}$  if there is an epimorphism
from $G_{1}$ to $G_{2}$, denoted $G_{1} \twoheadrightarrow G_{2}$.  This order is antisymmetric since each group in this set is Hopfian.  

Let $\mathcal{C} = \{G_{i}\}$ be a chain in $\Lambda.$  Each chain forms a direct system of groups over a totally ordered indexing set.  Any chain can be extended so that it starts with the element $F = G_{0}$.  Let $G_{\infty} = \varinjlim G_{i}$. 

By Lemma~\ref{RGlimit}, $RG_{p}(G_{\infty}) \geq \limsup RG_{p}(G_{i}) \geq \alpha$. Let $\widetilde{G}_{\infty}$ be the $p$-residualization of $G_{\infty}$.  By Lemma~\ref{residppropcomp}, $RG_{p}(\widetilde{G}_{\infty}) = RG_{p}(G_{\infty})$.  Therefore, $RG_{p}(\widetilde{G}_{\infty}) \geq \alpha$ and $\widetilde{G}_{\infty}$ is residually-$p$.  Moreover, for each $i$, $G_{i} \twoheadrightarrow G_{\infty}$ and in particular $F \twoheadrightarrow G_{\infty} \twoheadrightarrow \widetilde{G}_{\infty}$.  Thus $\widetilde{G}_{\infty} \in \Lambda$ and  $G_{i} \succcurlyeq \widetilde{G}_{\infty}$ for each $i$.  Thus, each chain $\mathcal{C}$ in $\Lambda$ has a lower bound in $\Lambda$ and therefore by Zorn's Lemma, $\Lambda$ has a minimal element, call it $\Gamma$.  

Since $\Gamma$ and its $p$-residualization $\widetilde{\Gamma}$ have the same $p$-gradient and $\Gamma$ surjects onto $\widetilde{\Gamma}$, it implies that $\widetilde{\Gamma} \in \Lambda$ and $\Gamma \succcurlyeq \widetilde{\Gamma}$.  Thus $\Gamma$ must be residually-$p$, otherwise $\widetilde{\Gamma}$ contradicts the minimality of $\Gamma$.  

\noindent \textbf{Note:}  $\Gamma$ does not have finite exponent. 

If $\Gamma$ had finite exponent then since $\Gamma$ is finitely generated and residually finite it must be finite by the positive solution to the Restricted Burnside Problem \cite{Zelmanov}.  This would imply $RG_{p}(\Gamma) <0$ by Corollary~\ref{finitepRG}.  This contradicts that $\Gamma$ is in $\Lambda$. 

Therefore, $\Gamma$ is a finitely generated residually-$p$ group with infinite exponent such that $RG_{p}(\Gamma) \geq \alpha$.  

\noindent \textbf{Claim:} $RG_{p}(\Gamma) = \alpha$. 

Assume not.  Then there exists a $k \in \mathbb{N}$ such that $RG_{p}(\Gamma) - \frac{1}{p^{k}} \geq \alpha$.  Since $\Gamma$ is residually-$p$, the order of every element is a power of $p$ and since $\Gamma$ has infinite exponent, there exists an $x \in \Gamma$ whose order is greater than $p^{k}$.

Consider $\Gamma' = \Gamma / \langle\langle x^{p^{k}} \rangle\rangle$.  Since $x^{p^{k}} \neq 1$ it implies that $\Gamma' \not\simeq \Gamma$.  By Theorem~\ref{pRGlowerbound}, $RG_{p}(\Gamma') \geq RG_{p}(\Gamma) - \frac{1}{p^{k}} \geq \alpha$.  If $\Gamma'$ is not residually-$p$, replace it with its $p$-residualization, which will have the same $p$-gradient. Then $\Gamma' \in \Lambda$ and $\Gamma \succcurlyeq \Gamma'$, which contradicts the minimality of $\Gamma$.  
\end{proof}

The result of Theorem~\ref{arbitraryRGp} can be strengthened without much effort.  

\begin{theorem}
\label{arbitrarytorsionRGp}
Fix a prime $p$.  For every real number $\alpha>0$ there exists a finitely generated residually-$p$ torsion group $\Gamma$ such that $RG_{p}(\Gamma) = \alpha$.
\end{theorem}

\begin{proof}
Barnea and Schlage-Puchta showed in Corollary 4 of \cite{BarneaPuchta}, that for any $\alpha >0$ there exists a torsion group $\mathcal{G}$ with $RG_{p}(\mathcal{G}) \geq \alpha$.
Applying the construction in Theorem~\ref{arbitraryRGp}, replacing the free group $F$ with the $p$-residualization of $\mathcal{G}$, will result in a group $\Gamma$ that is torsion, residually-$p$, and $RG_{p}(\Gamma) = \alpha$.  
\end{proof}

Y. Barnea and J.C. Schlage-Puchta \cite{BarneaPuchta} proved a result similar to Theorem~\ref{arbitrarytorsionRGp} (inequality instead of equality) albeit in a slightly different way.

\section{Applications}
The construction given in Theorem~\ref{arbitraryRGp} has a few immediate applications.  First, it is noted that Theorem~\ref{arbitrarytorsionRGp} gives a known counter example to the General Burnside Problem.  The second application is more general and shows that there exist uncountably many pairwise non-commensurable groups that are finitely generated, infinite, torsion, non-amenable, and residually-$p$. 

This paper has been concerned with the (absolute) rank gradient and $p$-gradient. There is, however, a related notion of rank gradient and $p$-gradient of a group with respect to a lattice of subgroups.  A set of subgroups, $\{H_{i}\}$, is called a \textit{lattice} if the intersection of any two subgroups in the set is also in the set.  In particular, a descending chain of subgroups is a lattice.

\begin{definition}  
\begin{enumerate}
\item The \textit{rank gradient relative to a lattice} $\{H_{i}\}$  of finite index subgroups is defined as  $\ds RG(G, \{H_{i}\}) = \inf_{i} \frac{d(H_{i}) - 1}{[G : H_{i}]}.$
\item The \textit{$p$-gradient relative to a lattice} $\{H_{i}\}$ of normal subgroups of $p$-power index is defined as $\ds RG_{p}(G, \{H_{i}\}) = \inf_{i} \frac{d_{p}(H_{i}) - 1}{[G : H_{i}]}.$
\end{enumerate}
\end{definition}

The following theorem was proved by Abert, Jaikin-Zapirain, and Nikolov in \cite{AJN}. Lackenby first proved the result for finitely presented groups in \cite{LacExpanders}.  

\begin{theorem}(Abert, Jaikin-Zapirain, Nikolov)
\label{RGamenable}
Finitely generated infinite amenable groups have rank gradient zero with respect to any normal chain with trivial intersection.
\end{theorem}

As a simple corollary, we provide a corresponding, albeit weaker, result concerning $p$-gradient.  

\begin{corollary}
\label{RGpamenable}
If $RG_{p}(G) >0$ for some prime $p$, then G is not amenable.  
\end{corollary}

\begin{proof}
Let $G$ be a finitely generated group with $RG_{p}(G)>0$.  Let $\widetilde{G}$ be the $p$-residualization of G.  Then $0<RG_{p}(G) = RG_{p}(\widetilde{G})$.  Let $\{H_{i}\}$ be a descending chain of normal subgroups of $p$-power index in $\widetilde{G}$ which intersect in the identity.  Then, 
\[ 0<RG_{p}(\widetilde{G}) \leq \inf_{i} \frac{d_{p}(H_{i})-1}{[\widetilde{G}:H_{i}]} \leq \inf_{i}\frac{d(H_{i})-1}{[\widetilde{G}:H_{i}]} = RG(\widetilde{G}, \{H_{i}\}). \]
Therefore, $\widetilde{G}$ is not amenable by Theorem~\ref{RGamenable}.  This implies that $G$ is not amenable since a quotient of an amenable group is amenable.
\end{proof}

The application of the construction used in Theorem~\ref{arbitraryRGp} concerning commensurable groups is given below.  

\begin{definition}
Two groups are called \textit{commensurable} if they have isomorphic subgroups of finite index.  
\end{definition}

The following lemma is straightforward.  

\begin{lemma}
Fix a prime $p$.  Let $G$ be a $p$-torsion group (every element has order a power of $p$).  Then every finite index subgroup $H \leq G$ is subnormal of $p$-power index.  
\end{lemma}

\begin{theorem}
There exist uncountably many pairwise non-commensurable groups that are finitely generated, infinite, torsion, non-amenable, and residually-p.  
\end{theorem}

\begin{proof}
Let $p$ be a fixed prime number.  By Theorem~\ref{arbitrarytorsionRGp} it is known that for every real number $\alpha>0$ there exists a finitely generated residually-$p$ infinite torsion group, $\Gamma$, such that $RG_{p}(\Gamma)=\alpha$.  By Corollary~\ref{RGpamenable}  these groups are all non-amenable.  Since each of these groups is residually-$p$ and torsion, they are all $p$-torsion.  Thus, every subgroup of finite index in these groups is subnormal of $p$-power index. 

By Theorem~\ref{pindexRGp} if any two of these groups are commensurable, then the $p$-gradient of each group is a rational multiple of the other.  Since there are uncountably many positive real numbers that are not rational multiples of each other, the result can be concluded.
\end{proof}

\bibliographystyle{amsplain}
\bibliography{rgcitjgrouptheory}

\providecommand{\bysame}{\leavevmode\hbox to3em{\hrulefill}\thinspace}
\providecommand{\MR}{\relax\ifhmode\unskip\space\fi MR }
\providecommand{\MRhref}[2]{%
  \href{http://www.ams.org/mathscinet-getitem?mr=#1}{#2}
}
\providecommand{\href}[2]{#2}
\begin{thebibliography}{10}

\bibitem{AJN}
M.~Abert, A.~Jaikin-Zapirain, and N.~Nikolov, \emph{The rank gradient from a
  combinatorial viewpoint}, Groups Geom. Dyn. \textbf{5} (2011), 213--230.

\bibitem{AtiyahMacDonald}
M.~Atiyah and I.~MacDonald, \emph{Introduction to commuative algebra}, Westview
  Press, 1969.

\bibitem{BarneaPuchta}
Y.~Barnea and J.C. Schlage-Puchta, \emph{On $p$-deficieny in groups},
  arxiv.org/abs/1106.3255v1 (2012).

\bibitem{DDMS}
J.D. Dixon, M.P.F du~Sautoy, A.~Mann, and D.~Segal, \emph{Analytic pro-p
  groups}, Cambridge University Press, 1991.

\bibitem{LacExpanders}
M.~Lackenby, \emph{Expanders, rank and graphs of groups}, Israel J. Math.
  \textbf{146} (2005), no.~1, 357--370.

\bibitem{NikolovSegal}
N.~Nikolov and D.~Segal, \emph{Finite index subgroups in profinite groups}, C.
  R. Math. Acad. Sci. Paris \textbf{337} (2003), no.~5, 303--308.

\bibitem{OsinRGTorsion}
D.~Osin, \emph{Rank gradient and torsion groups}, Bull. Lond. Math. Soc.
  (2010).

\bibitem{Pichot}
M.~Pichot, \emph{Semi-continuity of the first $\ell^{2}$-betti number on the
  space of finitely generated groups}, Comment. Math. Helv. \textbf{81} (2006),
  643--652.

\bibitem{Ribes}
L.~Ribes and P.~Zalesskii, \emph{Profinite groups}, 2 ed., A series of modern
  surveys in mathematics, vol.~40, Springer, 2010.

\bibitem{Puchta}
J-C. Schlage-Puchta, \emph{A $p$-group with positive rank gradient}, J. Group
  Theory \textbf{15} (2012), no.~2, 261--270.

\bibitem{wilson}
J.~S. Wilson, \emph{Profinite groups}, Oxford Science Publications, 1998.

\bibitem{Zelmanov}
E.~I. Zel'manov, \emph{Solution of the restricted burnside problem for groups
  of odd exponent}, Izv. Math. \textbf{54} (1990), no.~1, 42--59.

\end{thebibliography}

\end{document}